\documentclass[a4paper,12pt]{article}
\usepackage{amssymb}
\usepackage{amsthm}
\usepackage{amsxtra}
\usepackage{amsfonts}
\usepackage{xr}
\usepackage{color}
\usepackage{soul}
\usepackage{lscape}
\usepackage{colordvi}
\usepackage[normalem]{ulem}
\begin{document}

\newcommand{\red}{\textcolor{red}}
\newcommand{\meet}{\cap}
%\definecolor{green}{rgb}{0.1,.5,.1}
%\newcommand{\green}{\textcolor{green}}
\newcommand{\blue}{\textcolor{blue}}
\newcommand{\magenta}{\textcolor{magenta}}
%\definecolor{grey}{rgb}{0.5,0.5,0.5}
%\newcommand{\grey}{\textcolor{grey}}
\newcommand{\sect}{\section}

\newtheorem{thm}{Theorem}[section]
\newtheorem{cor}[thm]{Corollary}
\newtheorem{lem}[thm]{Lemma}
\newtheorem{prop}[thm]{Proposition}
\newtheorem{propconstr}[thm]{Proposition-Construction}

\theoremstyle{definition}
\newtheorem{para}[thm]{}
\newtheorem{ax}[thm]{Axiom}
\newtheorem{conj}[thm]{Conjecture}
\newtheorem{defn}[thm]{Definition}
\newtheorem{defns}[thm]{Definitions}
\newtheorem{notation}[thm]{Notation}
\newtheorem{rems}[thm]{Remarks}
\newtheorem{rem}[thm]{Remark}
\newtheorem{question}[thm]{Question}
\newtheorem{example}[thm]{Example}
\newtheorem{problem}[thm]{Problem}
\newtheorem{exercise}[thm]{Exercise}
\newtheorem{ex}[thm]{Exercise}

\newcommand{\mL}{{\mathcal L}}
\newcommand{\mC}{{\mathcal C}}
\newcommand{\spec}{{\rm Spec}}

\overfullrule=0pt
%\noitemsep

\newcommand{\si}{\sigma}
\newcommand{\prf}{\smallskip\noindent{\it        Proof}. }
\newcommand{\call}{{\mathcal L}}
\newcommand{\nat}{{\mathbb  N}}
\newcommand{\la}{\langle}
\newcommand{\ra}{\rangle}
\newcommand{\inv}{^{-1}}
\newcommand{\ld}{{\rm    ld}}
\newcommand{\trdeg}{{\rm tr.deg}}
\newcommand{\rest}{{\lower       .25     em      \hbox{$\vert$}}}
\newcommand{\union}{\cup}
\newcommand{\zee}{{\mathbb  Z}}
\newcommand{\conc}{^\frown}
\newcommand{\acls}{{\rm acl}_\si}
\newcommand{\acl}{{\rm acl}}
\newcommand{\cls}{{\rm cl}_\si}
\newcommand{\cals}{{\mathcal S}}
\newcommand{\calx}{{\mathcal X}}

\newcommand{\cl}{{\rm cl}}
\newcommand{\loc}{{\rm locus}}
\newcommand{\calg}{{\mathcal G}}
\newcommand{\calq}{{\mathcal Q}}
\newcommand{\calr}{{\mathcal R}}
\newcommand{\caly}{{\mathcal Y}}
\newcommand{\aff}{{\mathbb A}}
\newcommand{\cali}{{\mathcal I}}
\newcommand{\calu}{{\mathcal U}}
\newcommand{\frob}{{\rm Frob}}
\newcommand{\Frob}{{\rm Frob}}
\newcommand{\fix}{{\rm Fix}}
\newcommand{\Fix}{{\rm Fix}}
\newcommand{\per}{{\rm per}}
\newcommand{\dcl}{{\rm dcl}}
\newcommand{\calm}{{\mathcal M}}

\font\helpp=cmsy5
\newcommand{\semdp}
{\hbox{$\times\kern-.23em\lower-.1em\hbox{\helpp\char'152}$}\,}

\newcommand{\stab}{{\rm Stab}}
\newcommand{\qfcb}{\hbox{qf-Cb}}
\newcommand{\perf}{^{\rm perf}}
\newcommand{\sipm}{\si^{\pm 1}}
\newcommand{\SU}{{\rm SU}}
\newcommand{\evSU}{{\rm evSU}}
\newcommand{\vdim}{{\rm vdim}}
\newcommand{\vlabel}{\label}

\newcommand{\tp}{{\rm tp}}
\newcommand{\qftp}{{\rm qftp}}

\title{Revisiting virtual difference ideals}
\author{Zo\'e Chatzidakis\thanks{partially supported  by ValCoMo (ANR-13-BS01-0006)}
{ (CNRS
    - Universit\'e Paris Cit\'e (UMR 8576))}\and Ehud
Hrushovski\thanks{The research leading to these results has received
  funding from the European Research Council under the European Unions
  Seventh Framework Programme (FP7/2007- 2013)/ERC Grant Agreement
  No. 291111} (University of Oxford)} 

% }
% \date{}
%\centerline{\today}

\maketitle
\begin{abstract} In difference algebra, basic definable sets correspond to prime ideals that are invariant under a structural endomorphism.  
The  main idea of \cite{CHP} was that  periodic prime ideals enjoy better geometric properties than invariant ideals; and to understand
a definable set, it is helpful to enlarge it by relaxing invariance to periodicity, obtaining better geometric properties at the limit.
The limit in question was an intriguing but somewhat ephemeral setting called virtual ideals.  However a serious technical error was 
discovered by Tom Scanlon's UCB seminar.   In this text, we correct the problem via two different routes.   We replace the faulty lemma by a weaker one,
that still allows  recovering all results of \cite{CHP} for all virtual ideals.    
In addition, we introduce a  family of difference equations (``cumulative" equations) that we expect to be useful more generally.  Results in \cite{CHld} imply  that cumulative 
equations   suffice
to coordinatize all difference equation.  For cumulative equations, we show that virtual ideals reduce to globally periodic ideals, thus providing a proof of Zilber's trichotomy for difference equations using periodic ideals alone.  \end{abstract}

\section*{Introduction}

Boris Zilber developed a geometric description of $\aleph_1$-categorical theories, having a trichotomy at its heart.   It is based on the dimension theory of Morley (shown to take finite values by Baldwin), but gives information of a radically new kind than an abstract dimension theory.   Intuitively, a model of the theory is coordinatized by geometries that have either a graph-theoretic nature, or derive from linear algebra, or belong to algebraic geometry.    Though it is only the minimal definable sets that are described in this way, Zilber (and later others) demonstrated an overwhelming effect on the structure globally.  

 Zilber conjectured that there is no fourth option.  This turned out to be incorrect at the precise level of generality of $\aleph_1$-categorical structures.   But it was established 
with   additional hypotheses of a topological nature \cite{HZ}, and moreover proved to be 
 meaningful and indeed to capture the nature of structures  far beyond
 strong minimality.     Appropriate versions hold  for compact complex
 manifolds, for differentially closed and separably closed  fields, for
 strongly minimal sets  interpretable in algebraically closed fields of characteristic $0$ (\cite{Ca}); the latter closes in char. 0  a line opened more than thirty years ago by Eugenia Rabinovich, in her  Kemerovo PhD with Zilber.   The trichotomy 
 is also meaningful  for unstable theories: see \cite{PS} for the o-minimal case.     Many applications depend on the trichotomy, including 
 Zilber's gem \cite{zilber-jacobian}.      For difference equations, applications to diophantine geometry  include  \cite{HMM}, \cite{S},  \cite{CH2}.
 
 Thanks to Zilber's philosophy, when we made our first steps in the structure of difference equations in \cite{CH}, we knew in advance what it is that we should aim to prove.   The methods were informed by finite-rank stability and the nascent generalization to simplicity.   But they also 
 relied strongly on ramification divisors, and thus applied only in characteristic $0$.  Our approach in \cite{CHP} to the positive characteristic case thus had to be different.  

The trichotomy results of \cite{HZ} are valid for    stable structures with a finite dimension assigned to definable sets, satisfying  a `dimension theorem'
controlling dimensions of intersectinos.  Now the model companion  ACFA of the theory of difference fields is not stable, nor does the geometry of finite dimensional
sets satisfy the dimension theorem:  the intersection of two such sets may have unexpectedly low dimension.
For instance, the naive intersection of two surfaces in 3-space over the fixed field of the automorphism $\si$
could be two lines interchanged by $\si$; within the fixed field their intersection point would be the only solution.
Both of these pathologies   are ameliorated   as one relaxes $\si$ to $\si^m$
(going from the equation $\si(x) =F(x)$ to $\si^m(x) = F^{(m)}(x)$.)  At the limit, one has a {\em virtual structure},
defined and studied in \cite{CHP}; under appropriate conditions, this structure is stable and the dimension theorem is valid.  
Proving this uses basic ideas from topological dynamics  to obtain recurrent points, that may not be periodic; see Lemma \ref{lem08a} for example.
Using a generalization of the Zariski geometries of \cite{HZ}, one can then deduce the   trichotomy theorem.  The concrete form it takes here allows analyzing  any difference equation  via a tower of equations over fixed fields  and equations of locally modular type. \\[0.05in] 
In 2015, however,  Tom Scanlon's Berkeley seminar recognized a problem with a key technical lemma, 3.7.  
We show below how
to prove a somewhat weaker version of this lemma: where the wrong lemma~3.7 asserted
a unique component
through a point, the corrected version,   Proposition
\ref{lem37}, implies that   the number of such components 
  is finite, indeed at most the degree of the normalization of the relevant variety in the base.  All the main results of the paper remain valid with the same set of ideas, but considerable reorganization is required.  
One role of the present paper is to provide a lengthy  erratum, explaining in detail how this may be done.
Parts of this paper are thus technical and  need  to be read
in conjunction with \cite{CHP}. However section 2, which contains the main correction and in particular the key dimension theorem,   is self-contained   in the sense of
quoting some results from \cite{CHP} but not requiring entering into their
proofs.\\[0.05in]

At the same time, we take the opportunity to present a setting (`cumulative equations') in which the limit structure is equivalent to an ordinary structure,
in the sense that the associated algebraic object is an ordinary ring
with its periodic ideals, rather than an abstract limit of such rings as
in the virtual case. 
 Results of \cite{CHld} imply that
this setting, while not fully general, suffices to coordinatize all difference equations.  It may be of interest for other
applications, in particular the study of limit structures for more
equations that are not necessarily algebraic over SU-rank one. \\
We expect that a trichotomy theorem can be proved for Zariski geometries based on Robinson structures.  This has so far been worked
out only in special cases; the most general treatment is contained in
the unpublished PhD thesis of Elsner \cite{E}. Consequently the trichotomy follows from
the basic cumulative case alone, though this is not the case for some of the
other results: for finer
statements such as a description of the fields definable in the limit
structures, both in \cite{CHP} and here, we use additional features
of the specific structure.   \\[0.1in]

 Let $S$ be a difference ring, 
generated by a finitely generated   ring $R$.   
The main idea of \cite{CHP} was that as $n$ becomes more and more divisible, 
more $\si^n$-ideals appear, and their structures becomes progressively smoother.  
%the structure of prime $\si^n$- ideals simplifies as $n$ becomes more and more divisible.  
However there is also a countercurrent at work:  the difference subring $R_{\si^n}$ of $(S,\si^n)$ generated by $R$ may become smaller.  This double movement leads to technical complexity.    
If, however,   $\si(R)$ is contained in the ring generated by $R$ and $\si^n(R)$ for any $n$, this problem does not arise.
    It is this behaviour (slightly generalized to fraction fields)  that we refer to as  {\em cumulative}.  It turns out that
    cumulative difference equations still represent  all isogeny classes, and allow for considerable simplification.
     
%
%The algebraic object associated with a virtual structure is not a ring, but
%a descending sequence of rings; this is because a side-effect of moving from $\si$-ideals to $\si^m$-ideals
%is that a $\si$-ring 
%
 %
%
We are very grateful to Tom Scanlon, his Berkeley group, and especially Alex Kruckman for identifying 
the error; and to the anonymous  referee for the careful reading and suggestions that have considerably improved the text.
 \\[0.1in]
 
{\bf Plan of the paper}. In section 1 we mainly recall definitions and notations from
\cite{CHP}. Section 2 contains the proof of Proposition 2.6 of
\cite{CHP}, as well as some useful auxiliary results and remarks. The
cumulative case is done in the first half, the general case in the second
half. Sections 3 and 4 are devoted to rereading \cite{CHP} and making
the necessary changes and adaptations: Section 3 deals with sections 2
to 4 of \cite{CHP}, and section 4 with the remainder of the
paper.

\section{Setting, notation, basic definitions}

\para{\bf Setting and notation}. In
what follows, $K$ will be a sufficiently saturated existentially closed
difference field, containing an algebraically closed difference
subfield $k_0$, and $\Omega$ a $|K|^+$-saturated existentially closed
difference field containing $K$. We will always work inside $\Omega$. \\[0.05in]
If $L$ is a field, then $L^s$ % , $L^{\rm perf}$
and $L^{\rm alg}$ denote the
separable % , perfect
 and algebraic closure of the field $L$. \\[0.05in]
 {\bf Conventions}. 
Unless otherwise stated, all difference fields and rings will be
{\em inversive}, i.e., the endomorphism $\si$ is an automorphism;  
in other words we  take a difference ring to be a commutative ring with a $\mathbb{Z}$-action.
%{\bf onto}.  
Similarly, all
 difference ideals will be {\em reflexive}, i.e.: if $(R,\si)$ is a
difference ring, a $\si$-ideal of $R$ is an ideal $I$ such that
$\si(I)=\si ^{-1}(I)=I$. \\[0.05in]
If $k$ is a difference field, $X=(X_1,\ldots,X_n)$, then $k[X]_\si$ will denote the inversive difference
domain $k[\si^i(X_j)\mid i\in\zee, 1\leq j\leq n]$ and $k(X)_\si$ its
field of fractions. Similarly if $a$ is a tuple in $\Omega$: $k[a]_\si$
and $k(a)_\si$ denote the inversive difference subring and subfield of
$\Omega$ generated by $a$ over $k$. Similar notations for difference
rings. If $a$ is an $n$-tuple, then $I_\si(a/k)=\{f\in
k[X_1,\ldots,X_n]_\si\mid f(a)=0\}$. If  $k(a)_\si$ has
  finite transcendence degree over $k$, the {\em limit degree} of $a$
  over $k$, denoted $\ld(a/k)$ or $\ld_\si(a/k)$, is
  $\lim_{n\rightarrow\infty}[k(a,\ldots,\si^{n+1}(a)):k(a,\ldots,\si^n(a))]$. \\[0.05in]  
If $A$ is a subset of a difference
ring $S$, then $(A)_{\si^m}$ will denote the (reflexive) $\si^m$-ideal
of $S$
generated by $A$. If $A\subset \Omega$, then $\cl_\si(A)$ denotes the
perfect closure of the 
difference subfield of $\Omega$ generated by $A$, $\acls(A)$ the
(field-theoretic) algebraic closure of $\cl_\si(A)$, and $\dcl_\si(A)$
the model-theoretic definable closure of $A$. If $A$ is a subring of a difference ring $S$, then
$A_\si$ will denote the (inversive) difference subring of $S$ generated
by $A$. \\
Recall that $\acls(A)$ coincides with the model-theoretic algebraic closure $\acl(A)$, and that independence
(in the sense of the difference field $\Omega$) of $A$ and $B$ over a
subset $C$ coincides with the independence (in the sense of ACF) of $\acl(A)$ and $\acl(B)$
over $\acl(C)$. \\
If $m\geq 1$, then $\Omega[m]$ denotes the $\si^m$-difference field
$(\Omega,\si^m)$. The languages $\call$ and $\call[m]$ are the
languages $\{+,-,\cdot, 0,1,\si\}$ and $\{+,-,\cdot, 0,1,\si^m\}$. We
view $\call[m]$ as a sublanguage of $\call$, and $\Omega[m]$ as a
reduct of $\Omega$. Recall that $\Omega[m]$ is also an existentially
closed saturated difference field, by Corollary~1.12 of \cite{CH}. If $a$ is a tuple of $\Omega$ and $k$ a difference
subfield of $\Omega$, then $\qftp(a/k)$ denotes the quantifier-free type
of $a$ over $k$ in the difference field $\Omega$, and if $m\geq
1$, then $\qftp(a/k)[m]$ denotes the quantifier-free type of $a$ over $k$
in the difference field $\Omega[m]$. Similarly, if $q$ is a
quantifier-free type over $k$, then $q[m]$ denotes the set of
$\call(k)[m]$ quantifier-free formulas implied by $q$.

\bigskip\noindent{\bf \large Basic and semi-basic types}

\begin{defns}\vlabel{def1} We consider quantifier-free types $p$, $q$, \dots, over the
  algebraically closed difference field $k_0$, and integers $m,n\geq 1$. 
\begin{enumerate}

\item $q$ satisfies (ALG$m$) if whenever $a$ realises $q$, then
  $\si^m(a)\in k_0(a)^{\rm alg}$. 

\item The {\em eventual SU-rank} of $q$, $\evSU(q)$, is
$\lim_{m\rightarrow \infty}\SU(q[m!])$, where $\SU(q[m!])$ (the SU-rank
of $q[m!]$) is computed in
the $\si^{m!}$-difference field $\Omega[m!]$. For more details, see section 1  in \cite{CHP},
starting with 1.10. Notation: $\SU(a/k_0)[n]:=\SU(q[n])$, computed in
the $\si^n$-difference field $\Omega[n]$ ($n\geq 1$, $a$ realising $q$). If $D$ is  a
countable union of $k$-definable subsets of some cartesian power of
$\Omega$, then evSU$(D)=\sup \{\evSU(a/k)\mid a\in D\}$. 
\item $p\sim q$ if and only if for some $m\geq 1$, $p[m]=q[m]$. The
  $\sim$-equivalence class of $p$ is denoted by $[p]$ and is called a
  {\em virtual type}. 
\item $\calx_p(K)$ denotes the set of tuples in $K$ which realise $p[m]$
  for some $m\geq 1$.  Similarly for $\calx_p(\Omega)$. We denote by
  $X_p$ the underlying affine variety, i.e., the Zariski closure of
  $\calx_p(\Omega)$ in affine space. 
\item A {\em basic} type is a quantifier-free type $p$ over $k_0$, with
  evSU-rank $1$,  which satisfies (ALG$m$) for some $m$. Note that if $p$ is basic, so is
  $p[n]$ for every $n$. 
\item A {\em semi-basic type} is a quantifier-free type $q$ such that if
  $a$ realises $q$, then there are tuples $a_1,\ldots,a_n\in k_0(a)^{\rm
    alg}$
  which realise basic types over $k_0$, are algebraically independent over $k_0$, and are such that $a\in
  k_0(a_1,\ldots,a_n)^{\rm alg}$. 
\item The quantifier-free type $q$ is {\em cumulative} if for some (any)
  realisation $a$ of $q$ and every $m\geq 1$, $\si(a)\in
  k_0(a,\si^m(a))$. Note that this implies that $k_0(a)_\si=k_0(a)_{\si^m}$ for
  any $m\geq 1$, and that (ALG$m$) is equivalent to (ALG$1$). 
\end{enumerate}
\end{defns} 

\begin{rems}\vlabel{remld} Let $k$ be an inversive difference field. 
\begin{enumerate}
\item We will often use the following equivalences, for a tuple $a$:
\begin{enumerate} 
\item[(i)] $[k(a,\si(a)):k(a)]=\ld(a/k)$.
\item[(ii)] The fields $k(\si(a)\mid i\leq 0)$ and $k(\si^i(a)\mid i\geq
  0)$ are linearly disjoint over $k(a)$.
\item[(iii)] $I_\si(a/k)$ is the unique prime $\si$-ideal of $k[X]_\si$
  extending the prime ideal 
  $\{f(X,\si(X))\in k[X,\si(X)]\mid f(a,\si(a))=0\}$ of $k[X,\si(X)]$  ($|X|=|a|$). 
\end{enumerate}
Note that these equivalent conditions on the tuple $a$ in the difference
field $\Omega$ also imply the analogous conditions for the tuple $a$ in
the difference field $\Omega[m]$ for $m\geq 1$ (use (ii)).  
% $\ld_{\si^m}(a/k)=[k(a,\si^m(a)):k(a)]$ for any $m\geq 1$ (use (ii)). 
\item Let $P$ be a prime ideal of $k[X,\si(X)]$ ($X$ a tuple of
 variables) and assume that
  $\si(P\cap k[X])=P\cap k[\si(X)]$. Then $P$ extends to a prime
  $\si$-ideal of $k[X]_\si$. We will usually use it with the prime ideal 
  $\si\inv(P)$ of $k[\si\inv(X),X)]$. 
\end{enumerate}

\end{rems}

\prf All these are straightforward remarks;  see also section 1.3 of
\cite{CHld} for the equivalence of (1)(i) and (ii), and sections 5.6 and
5.2 of \cite{Co} for the remaining items.

\para\vlabel{coord0}\vlabel{def2}{\bf Coordinate  rings associated to
  quantifier-free types} (See also (3.5) and (3.6) in \cite{CHP}). Let
$q$ be a quantifier-free type over $k_0$, in the tuple $x$ of variables,
fix a realisation $a$ of $q$. The {\em pair $(R_q, R_{q,\si})$ of coordinate rings associated to
$q$} is defined as follows: Let $k_0(x)_\si$ be  the fraction field of
$k_0[X]_\si/I_\si(a/k_0)$, $k_0(x)$ its  subfield {generated by $x$
  over $k_0$}. Then we define
the ring $R_q:=k_0(x)\otimes_{k_0}K$ and the $\si^m$-difference ring $R_{q,\si^m}:=k_0(x)_{\si^m}\otimes_{k_0}K$
for $m\geq 1$. We often denote $R_q$ and $R_{q,\si^m}$ by $K\{x\}$
and $K\{x\}_{\si^m}$, and define in an analogous way the coordinate rings
$k_1\{x\}$ and $k_1\{x\}_{\si^m}$ if $k_1$ is a difference field
containing $k_0$. \\[0.1in]
Given  semi-basic types $q_1(x_1),\ldots,q_n(x_n)$, we  take
the tensor product over $K$ of their coordinate rings, and call them the
coordinate rings associated to $(q_1,\ldots,q_n)$. So, we  have
$$R_{(q_1,\ldots,q_n)}=K\{x_1\}\otimes_K\cdots \otimes_K K\{x_n\}, \ 
R_{(q_1,\ldots,q_n),\si^m}=K\{x_1\}_{\si^m}\otimes_K\cdots \otimes_K
K\{x_n\}_{\si^m}.$$ 
To a semi-basic type $q$, we  associate three new pairs of coordinate
rings as follows. Say $q$ is realised by a tuple
$a$, and $a_1,\ldots,a_n$ are as in the definition of semi-basic given
above. We let $p_i=\qftp(a_i/k_0)$, $r=\qftp(a_1,\ldots,a_n/k_0)$ and
$s=\qftp(a,a_1,\ldots,a_n/k_0)$. Then we define
$$R^1_q=R_{p_1}\otimes_K\cdots\otimes_K R_{p_n}, \quad
R^1_{q,\si^m}=R_{p_1,\si^m}\otimes_K\cdots\otimes_K R_{p_n,\si^m}$$
$$R^2_q=R_r, \ \ R^2_{q,\si^m}=R_{r,\si^m},\ \ R^3_q=R_s, \ \ R^3_{q,\si^m}=R_{s,\si^m}.$$  
These rings depend on the choice of the tuples $a_1,\ldots, a_n$, but we
may fix once and for all these tuples.  Note that then $R^1_q\subseteq
R^2_q\subseteq R^3_q\supseteq R_q$, and that $R_q^2$ is a localisation
of $R_q^1$, and $R_q^3$ is integral algebraic over $R_q^2$ and over
$R_q$. Similar statements hold for the associated difference rings. 
If $q$ is basic, we define $R^i_q=R_q$ and $R^i_{q,\si^m}=R_{q,\si^m}$.
We
extend the notation to the more general coordinate rings
$R_{(q_1,\ldots,q_n)}$. \\[0.05in]
We say that a coordinate ring $R_\si$ satisfies (ALG$m$) or is
cumulative, if the semi-basic types involved in the definition of $R_\si$
all satisfy (ALG$m$) or are cumulative.  

\para{\bf Convention}.  From now on, all quantifier-free types will
satisfy (ALG$m$) for some $m\geq 1$, so that all coordinate rings will
satisfy (ALG$m$).

\begin{defns}\vlabel{def3} Let $(R,R_\si)$ be a pair of coordinate rings, as defined
  above, and $S$ a ring. 
\begin{enumerate}
\item 
 Let $P$ be a prime ideal of a ring $S$. 
The {\em dimension} of $P$, denoted by $\dim(P)$, is the Krull dimension of the ring $S/P$.  If $I$
 is an ideal of $S$, the {\em dimension} of $I$, $\dim(I)$,  is $\sup \{\dim(P)\mid
 P\supseteq I,\ P\in \spec(S)\}$. If $S=R$, then $\dim(P)$ coincides with  $\trdeg_K{\rm
   Frac}(R/P)$. 
\item Let $P$ be a prime ideal of a coordinate ring $R_\si$. The {\em
    virtual dimension of $P$}, denoted $\vdim(P)$, is $\dim(P\cap R)$. If
    $R_\si$ satisfies (ALG$m$), it coincides with $\dim(P\cap
    R_{\si^m})$. Similarly, if $I$ is an ideal {of $R_\si$}, then $\vdim(I)=\dim(I\cap R)$.

\item A {\em virtual [perfect][prime] ideal} of $R_\si$  is a  [perfect\footnote{A $\si$-ideal $I$ of a
    difference ring $R$ is perfect if whenever $a^n\si(a)\in I$, then
    $a\in I$.}][prime] (reflexive) $\si^m$-ideal of $R_{\si^m}$
  for some $m\geq 1$. 
\item A {\em [perfect][prime] periodic ideal} of $R_\si$ is a
  [perfect][prime] $\si^m$-ideal $I$ of 
  $R_\si$ for some $m\geq 1$. A priori, not all virtual ideals extend to
  periodic ideals. 

\item Let $I$ be an  ideal of $R$. We say that $I$ is {\em pure of
    dimension $d$} if all minimal primes over $I$ have dimension
  $d$. Let $I$ be an ideal of $R_\si$. We say that $I$ is {\em virtually pure
  of dimension $d$} if $I\cap R$ is pure of dimension $d$. 

\item Let $I$ be a  virtual ideal of $R_\si=K\{x\}_\si$. Then
  $V(I)$ is the subset of $K^{|x|}$ defined by: $a\in V(I)$ if and only
  if for some 
  $m\geq 1$, for each $h \in I \meet R_{\si^m}$, 
  viewed as a $\si^m$-polynomial, we have 
 % $\si(a)=0$.
  $ h(a,\si^m(a),\cdots)=0$.   Thus $V(I)$ stands in bijection with $\union_m Hom_{\si^m}(R_{\si^m}/I,K)$,
  where $Hom_{\si^m}$ refers to ring homomorphisms commuting with $\si^m$.
  %Equivalently,  $(x-a)_{\si^m} \supseteq I\cap R_{\si^m}$. Note
  %$(x-a)_{\si^m}\cap R_{\si^m}\supseteq I\cap R_\si^m$. 
  
  Note
  that if  $R_\si=R_q$ for some quantifier-free  type $q$,
  then $V(0)$ is precisely $\calx_q(K)$. We call $\vdim(0)$ (i.e. the Krull dimension of $R$)  the {\em
  (virtual)  dimension of $q$}.    

\end{enumerate}
\end{defns}

\section{Existence theorems for periodic ideals}
The aim of this section is to give proofs of the results of \cite{CHP}
needed towards the proof of the trichotomy in positive
characteristic, and in particular the very important Proposition 2.6 of
\cite{CHP}. We try to follow the plan  % (exposition?)
of \cite{CHP}, and
will occasionally refer to it. While the results of chapter 2 are indeed
correct, the problem is that our coordinate rings do not satisfy the
required hypotheses. The mistake  appears in Lemma 3.7. \\[0.1in]

\ 

\noindent
{\bf \large Assumptions}\\
The
coordinate rings we consider are those associated to tensor products of
coordinate rings of semi-basic types whose corresponding basic types
have virtual 
dimension $e$, {\em for some fixed integer $e\geq 1$}. A typical pair of coordinate rings will be denoted $(R,R_\si)$,
without reference to the types involved in the construction.

As for types, we declare two virtual prime ideals $P,Q$ equivalent, and write  $P\sim
Q$, if   for some $m\geq 1$, $P\cap R_{\si^m}=Q\cap R_{\si^m}$.  We retain however definition 1.2(3) of virtual prime ideals; 
 the equivalence classes will be called virtual prime ideal classes.
 %(compare with} Definition 1.2(3)).
 
\begin{prop}\vlabel{prop24} {\rm (Addendum to Proposition 2.4 of
    \cite{CHP})} Let $(R,R_\si)$ be a pair of coordinate rings.
  \begin{enumerate}
\item  Let $P$ and $Q$ be virtual prime ideals. If $V(P)=V(Q)$, then $P\sim
Q$.   
\item Let $P$ be a prime $\si^m$-ideal of $R_{\si^m}$. Then for some
  $\ell>0$, $P$ extends to a  prime $\si^\ell$-ideal $Q$ of $R_\si$. In
  particular, since 
  $V(Q)=V(P)$, this shows that every    set defined by a virtual prime ideal is also
  defined by a periodic prime ideal  of $R_\si$;
  i.e. every prime periodic ideal of $R_{\si^m}$ is equivalent to a  prime periodic ideal of $R_{\si}$.

\end{enumerate}
\end{prop}

\prf 
(1) We may assume that $P$ and $Q$ are prime $\si$-ideals and that $R$
satisfies (ALG$1$). Choose a (small) subfield $k_1$ of $K$ such that for
any $m\geq 1$, $P\cap
R_{\si^m}$ and $Q\cap R_{\si^m}$ are  generated by their intersection
with $k_1\{x\}_{\si^m}$ $(x$ the variables of $R$). By saturation of $K$, it contains a point $a$
which is a generic point of $V(P)$ over $k_1$, i.e., with
$\trdeg(k_1(a)/k_1)=\dim(P)$. Then $a$ is in $V(Q)$, whence
$\dim(Q)\geq \dim(P)$, and the symmetric argument tells us that these
dimensions are equal, and that $a$ is a generic of $V(Q)$ over
$k_1$. Let $\ell$ be divisible by $m$ and such that $P\cap R_{\si^\ell}$ and $Q\cap
R_{\si^\ell}$ are prime $\si^\ell$-ideals contained in $(x-a)_{\si^\ell}$. Then 
$I_{\si^\ell}(a/k_1)=P\cap k_1\{x\}_{\si^\ell}=Q\cap k_1\{x\}_{\si^\ell}$, which
shows that $P\sim Q$. \\[0.05in]
(2) Let $\varphi: R_{\si^m}\to \Omega$ be a
$K$-homomorphism of $\si^m$-difference rings with kernel $P$. If
$p_1(x_1),\ldots,p_n(x_n)$ are the semi-basic types associated to
$R_\si$, then $R_\si=k_0(x_1)_\si\otimes_{k_0} \cdots \otimes
_{k_0}k_0(x_n)_\si\otimes_{k_0}K$, and $R_{\si^m}$ corresponds to the
subring $k_0(x_1)_{\si^m}\otimes_{k_0}\cdots \otimes
_{k_0}k_0(x_n)_{\si^m}\otimes_{k_0}K$. Our map $\varphi$ is entirely
determined by its restrictions to each of the factors of the tensor
product, and for $i=1,\ldots,n$, we let $\varphi_i$ denote the restriction of
$\varphi$ to $k_0(x_i)_{\si^m}$. Since $k_0(x)_\si$ is finitely
generated over $k_0(x)_{\si^m}$,  Proposition 1.12(3) of \cite{CHP}
gives that for some  $\ell>0$ divisible
by $m$, the $\si^\ell$-embeddings $\varphi_i:k_0(x_i)_{\si^m}\to \Omega$
extend to $\si^\ell$-embeddings $\psi_i:k_0(x)_\si\to \Omega$ for $i=1,\ldots,n$. Then
define $\psi=\psi_1\otimes \psi_2\otimes \cdots \otimes
\psi_n\otimes id_K$, and take $Q=\ker \psi$.

\begin{lem}\vlabel{lem0} Let $R_\si$ be a coordinate ring, and $S_\si=R[c]_\si$ a
  difference ring, with $S=R[c]$ integral algebraic (and finitely
  generated) over $R$.  If
  $P$ is a prime $\si$-ideal of $R_\si$, then for some $\ell\geq 1$,
  $P\cap R_{\si^\ell}$ extends to a 
  prime $\si^\ell$-ideal of $S_{\si^\ell}$. 
\end{lem}

\prf Replacing $\si$ by $\si^m$ for some $m$, we may assume that $R_\si$ satisfies (ALG$1$).\\[0.1in]
{\bf Claim}. There is $m\geq 1$ such that for any $\ell\geq 1$, if
$R'=R[\si(R),\ldots,\si^m(R)]$, then $P\cap R'_{\si^\ell}$ is the unique
prime $\si^\ell$-ideal of $R'_{\si^\ell}$ which extends 
$P\cap R'[\si^\ell(R')]$. \\[0.05in]
Indeed, let $a\in\Omega$ be such that ${\rm Frac}(R_\si/P)\simeq K(a)_\si$, and
choose $m$ such that
$[K(a,\ldots,\si^{m+1}(a)):K(a,\ldots,\si^{m}(a))]=\ld(a/K)$. Then if
$b=(a,\ldots,\si^m(a))$, we  have $\ld(b/K)=\ld(a/K)$ and for
$\ell\geq 1$, $\ld_{\si^\ell}(b/k_0)=[K(b,\si^\ell(b)):K(b)]$. \\
The claim now follows by {the equivalences given in} Remark \ref{remld}(1).\\[0.1in] 
For $n\geq 0$, let $R(n)$, resp. $S(n)$, denote the subring of $R_\si$,
$S_\si$ generated by $\si^i(R)$, 
$\si^i(S)$, $-n\leq i\leq n$. Then each $S(n)$ is Noetherian, integral
algebraic over $R(n)$,  
$S_\si=\bigcup_{n\in\nat} S(n)$, and we
have  a natural map $\spec (S_\si)\to
\prod_{n\in\nat}\spec(S(n))$. For each $n\in\nat$, the set $X_{n}$ of prime ideals of
$S(n)$ which extend $P\cap R(n)$ is finite and non-empty, and the
natural map $\spec (S(n+1))\to \spec( S(n))$ sends $X_{n+1}$ to $X_n$. Hence
$X:=\lim_\leftarrow X_{n}$ is a closed, compact, non-empty subset of
$\prod_{n\in\nat} X_{n}$, and is the set of prime ideals of $S_\si$
which extend $P$.  As each $X_{n}$ is finite,
and the set $X$ is stable under the (continuous) action of $\si$ on
$\spec (S_\si)$,  $X$ contains a recurrent point, $Q$. Let $m$ be given by
the claim, and consider $S(m)$. Then for some $\ell\geq 1$, we have
$\si^\ell(Q)\cap S(m)=Q\cap S(m)$, and therefore, using Remark \ref{remld}(2), there is a prime
$\si^\ell$-ideal $Q'$ of $S(m)_{\si^\ell}$ such that $$Q'\cap
S(m)[\si^{-\ell}(S(m))]=Q\cap S(m)[\si^{-\ell}(S(m))].$$ As $Q$ contains $P\cap
R'[\si^{-\ell}(R')]$ and has the same dimension, by the claim $Q'$ must
extend $P\cap R'_{\si^\ell}$, and therefore also $P\cap R_{\si^\ell}$. 

\begin{rem}\vlabel{evsu1} A consequence of our hypothesis on the
  dimension of the basic types is as follows: Let $P$ be a virtual prime ideal of
$R_\si$. Then $\dim(P\cap R)$ is divisible by $e$.  Indeed, choose $m$ such
that $P\cap R_{\si^m}$ is a prime $\si^m$-ideal of $R_{\si^m}$ and $R_\si$ satisfies
(ALG$m$). We may assume that $m=1$. We use the notation and definition
of \ref{def2}, and recall that $R^3$ is
finite  integral algebraic over $R$. Thus, by  Lemma \ref{lem0}, $P\cap R_{\si}$ extends to a periodic prime ideal of
$R^3_{\si}$. This means that 
 Frac$(R_{\si}/P\cap R_{\si})$ is equi-algebraic
over $K$ to 
a difference field  which is
generated over $K$ by realisations of basic types of 
dimension $e$. {Since basic types have evSU-rank $1$, these realisations
may be taken independent, and therefore}  $\trdeg_K({\rm Frac}(R_{\si}/P\cap R_{\si}))$
is a multiple of $e$, so that $\dim(P\cap R_{\si})$ is a multiple of
$e$. 
As $R_{\si}$ is integral algebraic over $R$, 
$\dim(P\cap R)$ is a multiple of $e$. 
\end{rem}

\bigskip\noindent
{\bf \large The basic cumulative case}\\[0.05in]
%
% \begin{lem}\vlabel{lem0b} Assume that $R_\si$ is a tensor product of
%   coordinate rings associated to cumulative quantifier-free types, and let
%   $S_\si$ be the quotient of $R_\si$ by some prime $\si$-ideal of $R_\si$. Then
%   $S_\si$ is also cumulative.
% \end{lem}
%
% \prf Write $R_\si=K\otimes k_0\{x_1\}_\si\otimes \cdots \otimes
% k_0\{x_n\}_\si$, and let $S_\si$ be as above. Noting that the natural map $R_\si\to S_\si$ restricts to
% isomorphisms on the difference subfields $k_0(x_i)_\si$ and $K$ gives
% the desired conclusion.\\[0.2in] 
We will now prove some results in the particular case when our
coordinate rings are tensor products of coordinate rings of {\bf basic 
  cumulative} types; this assumption holds until \ref{propdima}. The proof in the general case follows the same lines,
but is slightly more involved. \\
 
Note that the
assumptions  imply that all coordinate rings satisfy ALG$1$, that all
virtual ideals are periodic, and that $\sim$ coincides with equality. 

\begin{lem} \vlabel{lem05a} Let $I$ be an ideal of $R$ 
  of dimension 
  $d$. Then there are only finitely many periodic prime
  ideals  of
  $R_\si$ which contain 
$I$ and are of dimension $d$.
\end{lem}

\prf A prime ideal
of $R_\si$ which contains $I$ and is of  dimension $d$ must
extend a prime ideal $P$ of $R$ of dimension $d$ containing $I$. As $R$ is Noetherian,
there are only finitely many such  prime ideals, and 
we may therefore assume that $I=P$ is prime, and extends to a periodic
prime ideal of $R_\si$. \\[0.05in] 
Then Proposition~3.10 of \cite{CHP}, together
with Proposition \ref{prop24}, 
 gives 
the result.

\begin{cor}\vlabel{cor05a}  Let $I$ be an ideal of $R_\si$ of  dimension
  $d$. Then there are only finitely many periodic prime ideals of $R_\si$
  of  dimension $d$  containing $I$.
\end{cor}

\prf Such an ideal contains in particular $I\cap R$. The result follows
from Lemma \ref{lem05a}. 

\begin{cor}\vlabel{cor06a} Let $I$ be an ideal of
  $R_\si$ of  dimension $d$. Then there are periodic prime ideals $P_1,\ldots,P_s$ of
  $R_\si$ of  dimension $d$, and a finite subset $F$ of $I$, such
  that
if $P$ is a periodic prime ideal of $R_\si$ which contains $F$ and is of  dimension $d$, then $V(P)=V(P_i)$ for some $i$.
\end{cor}

\prf 
By Lemma \ref{lem05a}, if $F$ is a finite subset of $R_\si$ which
generates an ideal of dimension $d$ and $\per(F)$ denotes the set of
prime periodic ideals of $R_\si$ containing $F$ and of dimension $d$, then
$\per(F)$ is finite. Take a sufficiently large finite $F$ such that
$\per(F) =\per(I)$. 

\begin{lem} \vlabel{lem07a} Let $I$ be a periodic ideal of
  $R_\si$ of  dimension $d$. Then $I$ is contained in a periodic prime
  ideal of $R_\si$ of dimension $d$.\end{lem}

\prf We may assume that $I=\si(I)$. Let $F\subset I$ and $P_1,\ldots,P_s$
be given by
Corollary \ref{cor06a}. Let $X$ be the set of prime
ideals of $R_\si$ of dimension $d$  containing $I$, and for $n\in\nat$,
let $R(n)$ be the subring of $R_\si$ generated by $\si^i(R)$, $-n\leq
i\leq n$, and $X_n$ be the set  of prime ideals of $R(n)$ containing $I\cap R(n)$ and of dimension
$d$. Each $X_{n}$ is finite, non-empty, and we have  natural maps
$X_{n+1}\to X_n$. Hence, $X=\lim_\leftarrow X_n$ is non-empty and compact. The automorphism $\si$ acts continuously
on $X$, and therefore has a recurrent point $Q$. Let $n$
be such that $R(n)$ contains $F$. Then for some $m>0$,  we have $Q\cap R(n)=\si^m(Q)\cap
R(n)$. By Remark \ref{remld}(2),  there is a prime $\si^m$-ideal $Q'$ of $R(n)_{\si^m}$ which
extends  $Q\cap R(n)[\si^{-m}(R(n))]$. But $R(n)_{\si^m}=R_\si$, and
because $Q'$ contains $F$ and has dimension $d$, it must contain $I$. 

\begin{lem}\vlabel{lem08a} Let $I$ be a periodic
   ideal of $R_\si$, with $I\cap R$ 
  pure of dimension $d$. Then there are  periodic prime ideals
  $P_1,\ldots,P_s$ of virtual dimension $d$, such that $V(I)=V(P_1)\cup
  \cdots\cup V(P_s)$. 
\end{lem}

\prf We already know, by Lemma \ref{lem05a} (and Proposition \ref{prop24}), that $V(I)$ has only finitely
many irreducible components of dimension $d$, say $V(P_1),\ldots,V(P_s)$. It therefore suffices to
show that every point of $V(I)$ is in one of these components. Assume
this is not the case, let $a\in
V(I)$, and $m\geq 1$ such that $I$ is a 
$\si^m$-ideal and $Q=(x-a)_{\si^m}\supseteq I$. Without loss of generality, $m=1$. For
$n\in\nat$, let $R(n)$ be the subring of $R_\si$ generated by the rings
$\si^i(R)$, $-n\leq i\leq n$. Then for each $n\in\nat$, the  ideal $I\cap R(n)$ is pure of dimension $d$, and
therefore, the set $X_{n}$ of prime ideals $P$ of $R(n)$ of dimension $d$
containing $I\cap R(n)$ and contained in $Q$ is finite,
non-empty. Moreover, if
$P\in X_{n+1}$, then $P\cap R(n)\in
X_n$. Hence, the compact subset $X=\lim_{\leftarrow}X_{n}$ of
$\spec(R_\si)$ is non-empty. It is  the set of prime ideals
of $R_\si$ of dimension $d$, containing $I$ and contained in $Q$. 
Let $F$ be given by Lemma \ref{cor06a}, and $n$ such that $F\subset
R(n)$, and $Q$ does not contain any of the $P_i\cap R(n)$. As $\si$ acts
continuously on the compact set $X$, $X$ has a recurrent point, say
$P$. Then for some $m\geq 1$, $P\cap R(n)=\si^m(P)\cap R(n)$. As in the
proof of Lemma \ref{lem07a}, there is a prime $\si^m$-ideal $P'$ of
$R_\si$ which extends $P\cap R(n)[\si^{-m}(R(n)]$, and therefore has
dimension $d$, contains $I$ and is not in the finite set
$\{P_1,\ldots,P_s\}$. This gives us the desired contradiction. 

\bigskip

  We define a topology on $V$, taking the closed sets to be the sets $V(I)$.   (It is easy to see that the sets $V(I)$ are closed under intersections and under finite unions.)   
   Then 
when $s$ is taken minimal in Lemma \ref{lem08a}, the $V(P_i)$ are   the  {\em irreducible components} of $V(I)$.

\begin{lem}\vlabel{dim0} Write
  $R_\si=K\{x_1\}\otimes_K\cdots\otimes_KK\{x_m\}$, with $m\geq 2$, let
  $P$ be a prime $\si$-ideal of $R_\si$, and let $Q$ be the ideal 
   $Q=(x_1-x_2)_\si$  corresponding to the diagonal on $Spec K\{x_1\} \times Spec K\{x_2\}$, i.e. generated by the $x_{1,j} - x_{2,j}$.
  Then
  either $Q \subseteq P$, or every %non-empty 
  irreducible component of
  $V(P)\cap V(Q)$ has dimension $\dim(P)-e$.
\end{lem}

\prf Assume $Q \not\subseteq P$,    
and consider the $\si$-ideal $I=P+Q$.  Note that since $Q$ is generated by
elements of $R$, at least one of them is not in $P$; thus $I \meet R$ is strictly bigger than $P \meet R$;
so   each component of $I \meet R$ has dimension $<\dim(P)$. 

 Let
$R(n)$ {be} the subring of
$R_\si$
generated by $\si^i(R)$, $-n\leq i\leq n$ for $n\in\nat$. Then each $R(n)$ is the affine coordinate ring 
of a smooth variety, further localized (in fact in our construction of coordinate rings, {\em all} 
proper subvarieties over a certain field of definition were localized away; thus  including the 
singular locus of the variety.  (See the discussion given in (5.18) of \cite{CHP}).

 Hence the dimension theorem holds:  since $Q \cap R(n)$ has codimension $e$, 
 all minimal prime ideals   over
   $P\cap R(n)+Q\cap R(n)$ have dimension $\geq \dim (P)-e$.   
   
   Since $R$ is Noetherian, $I \meet R$ is finitely generated.  
Any finite set of elements of $I \meet R$ must already belong to  $P\cap R(n)+Q\cap R(n)$ for some $n$.  Since $R(n)$
is integral over $R$,  and the components of $P\cap R(n)+Q\cap R(n)$ have dimension $\geq  \dim
(P)-e$, it follows that every minimal prime of $I \meet R$ has dimension $\geq \dim
(P)-e$.  (The image of an irreducible variety under a morphism with finite fibers is an irreducible variety of the same dimension.)

 In particular, $I$ has  dimension $\delta \geq \dim(P)-e$.  By \ref{lem07a} some periodic prime ideal $P'$
containing $I$ has dimension $\delta$; by Remark \ref{evsu1},  $\delta$  as well as $\dim(P)$ must be a multiple of
$e$; we saw that $\delta < \dim(P)$, so  the only choice is  $\delta=\dim(P)-e$.    

Thus  $I \meet R$ is pure of dimension $\dim(P)-e$. 
  Hence \ref{lem08a} applies,
  and shows that  the components $V(P_1),\ldots,V(P_n)$ of $V(I)$ all have dimension exactly $d-e$.  
%
% Let $P_1,\ldots,P_s$ be as in 2.8; so the $V(P_i)$ are the components of $V(I)$.  
%We can find a finitely generated subring $R(n)$ of $R_\si$, such that the prime ideals $P_i \meet R(n)$
%are distinct, and no one contains another.   Then they are the minimal prime ideals containing $I \meet R(n)$.
%
% 
%
%Hence, every minimal
%prime ideal  %over
%  containing  $P+Q$ has dimension $\geq \dim(P)-e$. 
%  
%  
% By Lemma
%\ref{lem07a}, the $\si$-ideal $P+Q$ is contained in a prime periodic
%ideal $P'$ of 
%dimension $\dim(P+Q)$. By Remark \ref{evsu1}, \red{$\dim(P')$} must be a multiple of
%$e$, and this implies it must equal $\dim(P)-e$ {(because $\dim(P)$
%  is a multiple of $e$)}. {Moreover, no
%  minimal prime ideal of $(P\cap R)+(Q\cap R)$ of dimension
%  $>\dim(P)-e$  extends to a periodic ideal, and therefore all
%  irreducible components of $V(P)\cap V(Q)$ have dimension $\dim(P)-e$.} 

\begin{prop}\vlabel{propdima}   Let $P$ and $Q$ be periodic prime ideals
  of $R_\si$. Then every irreducible component  of $V(P)\cap V(Q)$ has
  dimension $\geq \dim(P)+\dim(Q)-\dim(0)$; it is determined by  a 
  periodic prime ideal of $R_\si$
 intersecting $R$ in  minimal prime ideals over $(P\cap R)+(Q\cap R)$.
   
\end{prop}
\prf  
 This can be deduced from  Lemma \ref{dim0} by reduction to an intersection with the diagonal $\Delta$ (identifying $V(P) \meet V(Q)$ with $P \times Q \meet \Delta$.)  
%
%follows immediately from Lemma \ref{dim0}, and the properties
%of dimension of ideals. 

\bigskip\noindent{\large \bf The general case}\\[0.05in]
The results in the cumulative case extend easily to the general case,
in most cases simply replacing equality of ideals by the equivalence
relation $\sim$. The fact that we consider also coordinate rings of
semi-basic types makes things a little more complicated, but Lemma
\ref{lem0} will be of use. 
Also, Proposition \ref{prop24} allows us to juggle between periodic and
virtual ideals. Recall our assumptions:\\[0.05in]
$(R,R_\si)$ is a tensor product of coordinate rings of semi-basic types,
and all associated basic types have virtual dimension $e$.

\begin{lem} \vlabel{lem05b} Let $I$ be an ideal of $R$, 
  of dimension 
  $d$. Then, up to $\sim$, there are only finitely many virtual prime
  ideals  of
  $R_\si$ which contain 
$I$ and are of virtual dimension $d$.
\end{lem}

\prf We may assume that $R_\si$ satisfies (ALG$1$). Then a prime ideal
of $R_\si$ which contains $I$ and is of virtual dimension $d$ must
extend a prime ideal $P$ of $R$ of dimension $d$ containing $I$. As $R$ is Noetherian,
there are only finitely many such  prime ideals, and 
we may therefore assume that $I=P$ is prime, and extends to a virtual 
prime ideal of $R_\si$. \\[0.05in] 
Let us first assume that the semi-basic types involved in $R_\si$ are all
basic. Then Proposition~3.10 of \cite{CHP}, together
with Proposition \ref{prop24}, 
 gives us
the result. \\[0.05in]
Let us now do the general case. We
will consider the rings $R^i$ introduced in \ref{def2}. Recall that
$R^1\subseteq R^2\subseteq R^3\supseteq R$. As $R^3_\si$ is integral
algebraic over $R_\si$, and satisfies (ALG$1$), Lemma
\ref{lem0} tells us that any virtual prime ideal of $R_\si$ extends to a virtual
prime ideal of $R^3_\si$. On the other hand, there are only finitely
many prime ideals of $R^3$ which extend $P$, so we may assume that
$R=R^3$, $R_\si=R^3_\si$. \\
The first case gives us that $P\cap R^1$
extends to finitely many prime virtual ideals of $R^1_\si$, up to
$\sim$, and by Proposition \ref{prop24}, we may assume they are
periodic. As $R^2$ and $R^2_\si$ are localizations of $R^1$ and $R^1_\si$
respectively, a periodic prime ideal of $R^1_\si$ extends to at most one
(periodic) prime ideal of $R^2_\si$. Say $Q$ is a
   prime
$\si^\ell$-ideal of $R^2_{\si^\ell}$ which extends $P\cap R^2$. Then
there are only finitely many prime ideals of $R^2_{\si^\ell}[R^3]$ which
extend $Q$, and by Lemma 3.9 of \cite{CHP}, to each of these
corresponds at most one (up to $\sim$) virtual
ideal of $R^3_\si$. Hence, up to $\sim$, there are
only finitely many virtual ideals of $R^3_\si$ extending $P$.

\begin{cor}\vlabel{cor05b}Let $I$ be an ideal of $R_\si$ of virtual dimension
  $d$. Then, up to $\sim$, there are only finitely many virtual prime ideals of $R_\si$
  of virtual dimension $d$ and which contain $I\cap R_{\si^m}$ for some $m>0$.
\end{cor}

\prf Such an ideal contains in particular $I\cap R$. The result follows
from Lemma \ref{lem05b}. 

\begin{cor}\vlabel{cor06b} Let $I$ be an ideal of
  $R_\si$ of virtual dimension $d$. Then there are periodic prime ideals $P_1,\ldots,P_s$ of
  $R_\si$ of virtual dimension $d$, and a finite subset $F$ of $I$, such
  that 
if $P$ is a periodic prime ideal  which contains $F$ and is of virtual dimension $d$, then $V(P)=V(P_i)$ for some $i$.
\end{cor}

\prf 
By \ref{lem05b}, if $F$ is a finite subset of $R_\si$ which
generates an ideal of dimension $d$ and $\per(F)$ denotes the set of
prime periodic ideals of $R_\si$ containing $F$ and of dimension $d$, then
$\per(F)/\!\!\sim$ is finite. Take a sufficiently large finite $F$ such that
$\per(F)/\!\!\sim =\per(I)/\!\!\sim$. 

\para{\bf Warning}. This set $F$ is not necessarily contained in $R$,
nor in $\bigcap_m R_{\si^m}$, unless $R_\si$ is cumulative. 

We will need a version of Lemma \ref{lem08a} without the purity assumption.   
We  claim  a weaker conclusion, 
namely that $V(I)$ is contained in some $V(P_i)$ of maximal dimension.

\begin{lem} \vlabel{lem07b} Let $I$ be a virtual ideal of
  $R_\si$ of virtual dimension $d$. Then there are $m\geq 1$ and a prime
  $\si^m$-ideal of $R_{\si^m}$ of dimension $d$ which contains $I\cap R_{\si^m}$.\end{lem}

\prf We may assume that $I=\si(I)$, and that
$R_\si$ satisfies (ALG$1$). Let $F\subset I$
be given by
Corollary \ref{cor06b}. Let $X$ be the set of prime
ideals of $R_\si$ of dimension $d$  containing $I$, and for $n\in\nat$,
let $R(n)$ be the subring of $R_\si$ generated by $\si^i(R)$, $-n\leq
i\leq n$, and $X_n$ be the set  of prime ideals of $R(n)$ containing $I\cap R(n)$ and of dimension
$d$. Each $X_{n}$ is finite, non-empty, and we have  natural maps $X\to
\prod_{n\in\nat}X_{n}$ and $X_{n+1}\to X_n$. The automorphism $\si$ acts continuously
on the compact set $X$, and therefore has a recurrent point $Q$. Let $n$
be such that $R(n)$ contains $F$. Then for some $m>0$,  we have $Q\cap R(n)=\si^m(Q)\cap
R(n)$. By Remark~\ref{remld}(2),  there is a prime $\si^m$-ideal $Q'$ of $R(n)_{\si^m}$ which
extends  $Q\cap R(n)[\si^{-m}(R(n))]$. Applying Proposition \ref{prop24}
to $R(n)_{\si^m}$, we obtain a prime $\si^\ell$-ideal $Q''$ of $R_\si$ which
extends $Q'$; then $Q''$ contains $F$ and has dimension $d$.

\begin{lem}\vlabel{lem37} {\rm (Correct version of Lemma 3.7 in \cite{CHP})}\,  Let $R$ be a domain which is
  integrally closed, let $k$ be a subfield of $R$, and $k_1$ an
  algebraic extension of $k$, and let $S=k_1\otimes_k R$. Let $Q$ be a
  prime ideal of $S$. 
\begin{enumerate}
\item There is a unique prime ideal of $S$ which intersect
  $R$ in $(0)$ and is contained in $Q$. 
\item If  $P'$ is a prime ideal
  of $S$ which intersects $R$ in $(0)$ and if $k_1$ is  separably
  algebraic over $k$,  then $S/P'$ is integrally closed. 
\end{enumerate}
\end{lem}

\prf 
For both (1) and (2), we may assume that $S$ is finitely generated over
$R$, i.e., that
$k_1$ is a finite extension of
$k$. Furthermore, observe that if $b\in S$, then 
$b^{p^n}$ belongs to the subring $ (k_1\cap k^s)\otimes_k R$ of $S$ for some $n$, and
that a prime ideal $P$ of $S$ contains $b$ if and only if its
intersection with $ (k_1\cap k^s)\otimes_k R$ contains $b^{p^n}$. I.e.,
the restriction map $ \spec(S)\to \spec((k_1\cap k^s)\otimes_k R)$ is a bijection. We may
therefore assume that $k_1$ is separably algebraic over $k$, of the form
$k[a]$ for some $a\in k_1$. \\
Let $f(T)$ be the
minimal monic polynomial of $a$ over $k$ and consider its
factorization $\prod_{i=1}^m g_i(T)$ over ${\rm Frac}(R)$ into monic
irreducible 
polynomials. Because $R$ is integrally closed, all $g_i(T)$ are in
$R[T]$ (see e.g. Thm~4, Ch~V \S3 in \cite{ZS}). Moreover, since $f$ is
separable, their coefficients are 
actually in the subfield $R\cap k^s$ of $R$, and if $i\neq j$, then $(g_i(T),g_j(T))=(1)$. Thus any prime ideal of
$S$, and in particular $Q$, contains one and  only one
of the elements $g_i(a)$, and the ideal of $S$ generated by $g_i(a)$ is 
 prime. (For this last assertion, use the fact that $g_i(T)$ is
 irreducible over ${\rm Frac}(R)$, and that $S\simeq R[T]/f(T)$). This
 shows (1).\\[0.05in] 
(2) Viewing $R$ as the coordinate ring of an affine  variety $V$ over
$k$, we know that $V$ is normal. A minimal prime ideal of $S$
 corresponds 
therefore to an irreducible component of the (non-irreducible) variety
$V_{k_1}$, and as the property of normality is a local property, each
component of $V_{k_1}$ is normal, i.e., with $P'$ as above, $S/P'$ is
integrally closed. Here we are using the fact that $k_1/k$ is separable,
so that the map $\spec( k_1)\to \spec( k)$ is \'etale and if $k_1/k$ is
finite, then  $S$ is a product of domains. \\
The fact that $R$ is not necessarily finitely generated over $K$ is
not important: it is a union of finitely generated $K$-algebras which
are integrally closed. 

\begin{prop}\vlabel{prop8} Let $(R,R_\si)$ be a pair of coordinate rings associated to
  semi-basic types satisfying (ALG$1$). Then $(R,R_\si)$ satisfies the following:  if
  $Q$ is a prime ideal of $R_\si$ and if $P$ is a prime ideal of
  $R$ which is contained in $Q\cap R$, then there are only finitely many
  prime ideals of $R_{\si}$ which extend $P$ and are contained in
  $Q$. \end{prop}

\prf Let $Q\subset R_{\si}=S$ be a prime ideal, let $P$ be
a prime ideal of $R$ such that $P\subseteq Q\cap R$. Let us first assume
that $R/P$ is integrally closed. Let $(x_1,\ldots,x_n)$ be the
coordinates corresponding to $R$, i.e.,
$R=K\{x_1\}\otimes_K\cdots\otimes_KK\{x_n\}$ and $K\{x_i\}=k_0(x_i)\otimes_{k_0}K$. Then
$$S=\bigl(\cdots\bigl(\bigl(R\otimes_{K\{x_1\}}K\{x_1\}_\si\bigr)\otimes_{K\{x_2\}}K\{x_2\}_\si\bigr)\cdots
\otimes_{K\{x_n\}}K\{x_n\}_\si\bigr).$$ We
know that each $K\{x_i\}_\si$ is integral algebraic over $K\{x_i\}$ (by (ALG$1$)). However, it
may not be separably integral algebraic. So, we will consider instead the ring 
$$S'= \bigl(\cdots\bigl(R\otimes_{K\{x_1\}}(K\{x_1\}_\si\cap K\{x_1\}^s)\bigr)\otimes_{K\{x_2\}}\cdots
\otimes_{K\{x_n\}}(K\{x_n\}_\si\cap K\{x_n\}^s)\bigr).$$
If $b\in S$, some $p^m$-th power of $b$ lies in $S'$, so that any
prime ideal of $S'$ extends uniquely to a prime ideal of $S$. It
therefore suffices to prove the result for $S'$. \\
Applying Lemma \ref{lem37} to
$k=K\{x_1\}$ and $S_1=R\otimes_{K\{x_1\}}(K\{x_1\}_\si\cap K\{x_1\}^s)$, we obtain that
  there is a unique prime ideal $P_1$ of $S_1$ which extends $P$ and is
  contained in $Q\cap S_1$. Furthermore, $S_1/P_1$ is integrally
  closed. Iterate the reasoning to obtain that there is a unique
  prime ideal $P_n$ of $S'$ which extends $P$ and is contained in $Q$
  (and furthermore, $S'/P_n$ is integrally closed).\\[0.05in]
In the general case, 
let $A$ be the integral closure of $R/P$. Because
$R/P$ is a localization of a finitely generated $K$-algebra, it follows that $A$
is a finite $R/P$-module  (see \cite{ZS}, Ch~V, \S4 Thm~9; observe also
that a localization of an integrally closed domain is integrally closed), and is
integral algebraic over $R/P$. So the map $\spec (A)\to \spec(
R/P)$ is finite, with fibers of size at most $g$ for some $g$. 
Hence, the prime ideal $Q/PS$ of $S/PS$ has exactly $s$  
extensions $Q_1,\ldots,Q_s$ to $\tilde S=(S/PS)\otimes_{R/P}A$, for some $s$ with $1\leq s\leq
g$. Let $P'$ be a prime ideal of $S$ extending $P$ and contained in $Q$; then $P'$ contains
$PS$, and therefore $P'/PS$ extends to a prime ideal $Q'$ of $\tilde S$; this
$Q'$  must be  contained in one of the $Q_i$'s. By the first case, this
determines $Q'$ uniquely, and therefore also $P'$. Hence $P$ has at most
$s$ extensions to prime ideals of $R_\si$ which are contained in $Q$. 
  
\begin{lem}\vlabel{lem08b} Let $I$ be a  virtual
  perfect ideal of $R_\si$, with $I\cap R$ 
  pure of dimension $d$. Then there are periodic prime ideals
  $P_1,\ldots,P_s$ of virtual dimension $d$, such that $V(I)=V(P_1)\cup
  \cdots\cup V(P_s)$. 
\end{lem}

\prf We already know, by Lemma \ref{lem05b}, that $V(I)$ has only finitely
many irreducible components of dimension $d$. It therefore suffices to
show that every point of $V(I)$ is in one of these components. Let $a\in
V(I)$, and $m\geq 1$ such that $R_\si$ satisfies (ALG$m$), $I\cap R_{\si^m}$ is a perfect
$\si^m$-ideal and $Q=(x-a)_{\si^m}\supseteq I\cap R_{\si^m}$. We
will work in $R_{\si^m}$, so without loss of generality, $m=1$. For
$n\in\nat$, let $R(n)$ be the subring of $R_\si$ generated by the rings
$\si^i(R)$, $-n\leq i\leq n$. Then for each $n\in\nat$, the  ideal $I\cap R(n)$ is pure of dimension $d$, and
therefore, the set $X_{n}$ of prime ideals $P$ of $R(n)$ of dimension $d$
containing $I\cap R(n)$ and contained in $Q$ is finite,
non-empty. Moreover, if
$P\in X_{n+1}$, then $P\cap R(n)\in
X_n$. Hence, the compact subset $X=\lim_{\leftarrow}X_{n}$ of
$\spec(R_\si)$ is non-empty. It is  the set of prime ideals
of $R_\si$ of dimension $d$, containing $I$ and contained in $Q$. If
$P\in X$, then $P\cap R$ belongs to the finite set $X_0$; hence, by
Lemma \ref{lem37}, $X$ is finite. On the other hand, $X$ is stable under
the (continuous) action of $\si$, because $I$ and $Q$ are $\si$-ideals. Hence, for
some $\ell$, $\si^\ell$ is the identity on $X$, i.e., all ideals in $X$
are prime $\si^\ell$-ideals. 

\begin{prop} \vlabel{prop26} {\rm (Proposition 2.6 in \cite{CHP})\, } Let $(R,R_\si)\in\calr$ be a pair of
  coordinate rings, and let
$P_1,P_2$ be two virtual prime ideals of $R_\si$.  Then $V(P_1)\cap
V(P_2)=V(I)$ for some virtual perfect ideal $I$. The irreducible components of 
$V(P_1)\cap V(P_2)$ correspond to virtual prime
ideals $Q_i$ with $Q_i\cap R$ minimal prime containing $P_1\cap R +P_2\cap R$. 
\end{prop}

\prf We may assume that $R_\si$ satisfies (ALG$1$), and that $P_1$ and
$P_2$ are prime $\si$-ideals. (In fact, at every stage of the proof, we
will allow ourselves to replace $R_\si$ by $R_{\si^m}$ so that our ideals
remain $\si$-ideals, and without explicitly saying so). For the first
assertion, it suffices to show that $V(P_1)\cap V(P_2)$ has only
finitely many irreducible components: if these are of the form $V(Q_i)$,
$i=1,\ldots,s$, for $Q_i$ a prime $\si^m$-ideal of $R_{\si^m}$, then one
takes $I=\bigcap_{i=1}^s Q_i$, a perfect $\si^m$-ideal of
$R_{\si^m}$ (which contains $P_1\cap R_{\si^m}+P_2\cap R_{\si^m}$).\\[0.05in]
If $V(P_1)\cap V(P_2)=\emptyset$, there is nothing to prove, so we will
assume  it is non-empty. The elements of $V(P_1)\cap V(P_2)$ are in
correspondence with the elements of $(V(P_1)\times V(P_2))\cap \Delta$,
where the corresponding pair of coordinate rings is
$(R_\si\otimes_KR_\sigma, R\otimes_KR)$, and $\Delta$ denotes the
diagonal of the underlying ambient set $V(0)\times V(0)$. The same
observation holds at the level of the Zariski closures. We will
therefore replace $P_1$ by the ideal $P$ of $R_\si\otimes_KR_\si$
generated by $P_1\otimes 1+1\otimes P_2$, and $P_2$ by the ideal
corresponding to $\Delta$, i.e., the ideal $I(\Delta)$ of $R_\si\otimes_KR_\si$
generated by all $a\otimes 1-1\otimes a$, for $a\in R_\si$. Write
$R_\si$ as the tensor product over $K$ of the rings $K\{x_i\}_\si$,
$i=1,\ldots,n$, with
$K\{x_i\}$ associated to the semi-basic type $q_i$. Then $\Delta=\bigcap
\Delta_i$, where $\Delta_i\subset V(0)\times V(0)$ is defined by
$x_i=x'_i$ inside $$S_\si=(K\{x_1\}_\si\otimes_K\cdots
\otimes_KK\{x_n\}_\si)\otimes_K(K\{x'_1\}_\si\otimes_K\cdots
\otimes_KK\{x'_n\}_\si).$$
It then suffices to show the result for $P+I(\Delta_1)$, then for each
$P'+I(\Delta_2)$ where $P'$ is a  prime periodic ideal minimal
containing $P+I(\Delta_1)$, etc. \\[0.05in]
Let us first assume that  $q_i$ is basic and that   $P$ does not contain
$I(\Delta_i)$. The proof is  very similar to the proof of  Lemma
\ref{dim0}, with small changes. Let $S=R\otimes_KR$,  
$S_\si=R_\si\otimes_K R_\si$, {and $S(n)\subset S_\si$ the subring} 
generated by $\si^i(S)$, $-n\leq i\leq n$ for $n\in\nat$. Reasoning as
in the proof of \ref{dim0}, {all minimal prime ideals over
$P+I(\Delta_i)$ have dimension $\geq \dim(P)-e$.} 
By Lemma
\ref{lem07b}, $P+I(\Delta_i)$ is contained in a prime periodic ideal $P'$ of
dimension $\dim(P+I(\Delta_i))$. By Remark \ref{evsu1}, $\dim(P+I(\Delta_i))$ must be a multiple of
$e$, and this implies it must equal $\dim(P)-e$. Hence all irreducible
components of  $V(P+I(\Delta_i))$ have
dimension $\dim(P)-e$. \\
 Note that the minimal virtual prime  ideals containing
$P+I(\Delta_i)$ do indeed extend minimal prime ideals over $P\cap
S+I(\Delta_i)\cap S$, since they have the same dimension. \\[0.05in]
We will now do the general case. As $R^3_\si$ is integral algebraic over
$R_\si$, we may assume that $R_{q_i}=R^3_{q_i}$, $R_{q_i,\si}=R^3_{q_i,\si}$, by Lemma
\ref{lem0}. Write the variables of $q_i$ as 
$(y,y_1,\ldots,y_r)$. Then $I(\Delta_i)$ is the intersection of the $r$ 
$\si$-ideals $(y_1-y'_1)_\si$, $(y_2-y'_2)_\si$, \dots,
$(y_r-y'_r, y-y')_\si$. The first $r-1$ of these ideals have dimension
$\trdeg_K(S)-e$ in $S_\si$; for the last one, work inside $S_\si/(y_1-y'_1,
y_2-y'_2, \dots,y_{r-1}-y'_{r-1})_\si$: then the minimal prime $\si$-ideals over $I(\Delta_i)/(y_1-y'_1,
y_2-y'_2)_\si, \dots,y_{r-1}-y'_{r-1})_\si$ all have dimension
$\trdeg_K(R_\si)$.  Apply the first case to these ideals to conclude.

\begin{cor}\vlabel{dim} {\rm (The dimension theorem - see 4.16 in \cite{CHP})} Let
  $P_1$ and $P_2$ be virtual prime ideals of $R_\si$, and let $n$ be the
  evSU-rank of $V(0)$. (I.e., there are
  exactly $n$ basic types which are associated to $R_\si$). Then all
  non-empty irreducible components of $V(P_1)\cap V(P_2)$ have evSU-rank
  $\geq (\dim(P_1)+\dim(P_2))/e -n$. 
\end{cor}

\section{Going through sections 2, 3 and 4 of \cite{CHP}}
We will describe which of the results of these three sections remain
true without changes, which ones are false or unnecessary, and which ones
need to be repaired. Note that while our coordinate rings are not
``friendly'' (because they do not satisfy $(*1)$), the assumption we
make on the semi-basic types considered are usually slightly stronger
than those made in 
the paper.  Unless preceded by ``the present'', references are to results in
\cite{CHP}. \\[0.1in]
{\bf \large Section 2}\\
We gave up on the idea of finding a
general setting (a modified version of friendliness satisfied by our
coordinate rings) in which one would
be able to prove the dichotomy 
theorem, and so in all the results, the hypotheses of friendliness
should be replaced by our hypotheses on semi-basic types: the associated
basic types all have dimension $e$. \\[0.05in]
Notation and definitions are given in more details in paragraphs 2.1
and 2.2, as well as some examples. Proposition 2.4 states the basic
results on the 
duality between sets $V(I)$ and virtual ideals. \\
Proposition 2.6 is the present  Proposition \ref{prop26}. The proof of
Proposition 2.8 goes through verbatim.\\[0.1in]
{\bf \large Section 3}.\\
Paragraphs (3.1) to (3.6) are definitions and notations. \\
Lemma 3.7 is {\bf false}, the correct version is given by the present  Lemma
\ref{lem37}(1), but is not enough to prove $(*1)$ for our coordinate
rings. Thus Proposition 3.8 is  false as well. \\
However, the proofs of Lemma 3.9 and Proposition 3.10 go through,
without change (except for a typo on line 4 of the proof of 3.10, it
should be $Q\cap
K[x_1,\ldots,x_r]_\si$). \\
Theorem 3.11 is implied by the present  Corollary \ref{cor05b}. \\
Proposition 3.12 goes through verbatim (note that the claim is the present 
Remark \ref{evsu1}). Note also that once more, Proposition 2.6 (i.e.,
the present  Proposition \ref{prop26}) is
instrumental. \\[0.1in]
{\bf \large Section 4}\\
Paragraph 4.1 consists of definitions and notations. \\[0.05in]
Proposition 4.2 remains true, but the proof needs to be slightly
modified (as it appeals to the false Lemma 3.7) towards the end. The
modification is as follows: we are in the situation of $R_\si$ satisfying
(ALG$1$), have chosen $a_1,\ldots,a_n,a\in V(P)$ such that the field of
definition of the ideal $P\cap R$ is
contained in $k_0(a_1,\ldots,a_n)$, and $a$ is generic over
$k_0(a_1,\ldots,a_n)$. By (ALG$1$) and the way our coordinate rings are
defined, we know that the ideal $I$ of $R_\si$ generated by $P\cap R$ is pure
of dimension $\dim(P)$. As $V(I)$ has finitely many irreducible
components and by genericity of $a$, $a$ is in only one irreducible
component of $V(I)$, and that component must be $V(P)$. Hence, for any
$\ell$, $P\cap R_{\si^\ell}$ is defined over
$\cl_{\si^\ell}(k_0,a,a_1,\ldots,a_n)$. \\[0.05in]
Corollary 4.2, Propositions 4.3, 4.4 and 4.5 go through without
change. (In the proof of 4.3, replace $(*1)$ by the present  Proposition
\ref{prop8})\\[0.05in] 
In 4.6, we will slightly strengthen the requirements and only consider
$0$-closed sets defined by {\em virtual 
  perfect} ideals. This is to ensure that they have only finitely many
irreducible components. \\[0.05in]
Proposition 4.7 remains true, with a slight change at the end of the
proof, similar to the one given for 4.3.\\
Proposition 4.9 and Lemma 4.10 go through without change. Note the
following consequences of Lemma 4.10, which are quite useful and were
not  stressed enough in the paper \cite{CHP}:\\[0.05in]
{\bf Corollaries of Lemma 4.10 of \cite{CHP}}. (1) {\em Let $d_1$ and $d_2$ be
  tuples of 
  realisations of basic types among $\{p_1,\ldots,p_n\}$. Then
  $\acl(d_1)\cap \acl(d_2)=\acl(e)$, where $e$ consists of
  realisations of types in $\{p_1,\ldots,p_n\}$.}\\
(2) {\em Let $b$ realise a tuple of semi-basic types, and $a\in \acl(b)$ be
such that $\qftp(a/k_0)$ satisfies (ALG$m$) for some $m$. Then
$\qftp(a/k_0)$ is semi-basic.}\\[0.05in]

\prf (2) Indeed, without loss of generality $b$
consists of realisations of basic types; take $b'$ realising $\qftp(b/a)$
and independent from $b$ over $a$. Then $a=\acl(b)\cap \acl(b')$ and we
may apply (1). \\[0.2in]
Let us now discuss Theorem 4.11. The set $\caly$  needs to be modified in the following
manner: \\
Condition (i) (stays the
same): for any semi-basic type $q$, $\calx_q(K)\subset \caly(K)$ or
$\calx_q(K)\cap \caly(K)=\emptyset$;\\
Condition
(ii) becomes: if 
$b\in\caly(K)^n$ for some $n$, and $a\in \acl(k_0b)$ is such that
$q=\qftp(a/k_0)$ satisfies (ALG$m$) for some $m$, then
$\calx_q(K)\subset \caly(K)$.  [The set $\caly$ was in fact incorrectly
defined in \cite{CHP}, and the current definition is the one which is  used in the
proof]. In the cumulative
case, we furthermore impose that all our semi-basic types are
cumulative.\\[0.05in]
Once this change done, the proof goes through, however one needs to pay
attention to a clash of notation: the tuple $d$ which appears on line 13 of
page 283 has nothing to do with the one discussed earlier in the proof;
it consists of realisations of basic types, and is independent from $c$
over $k_0$.   \\[0.05in] 
Proposition 4.12 of \cite{CHP} goes through verbatim, as well as Remark 4.14, Proposition
4.15 and the verification of the axioms for Zariski geometries given in
(4.16), for the set $\caly_b(K)=\bigcup_{p {\rm \ basic\ of\ dimension\
  }e}\calx_p(K)$. Note
that the present  Corollary \ref{dim} gives
us Corollary 4.16 of \cite{CHP} for semi-basic types. 

\section{Using the Zariski geometry to get the trichotomy}
The first paragraphs of chapter 5 of \cite{CHP} introduce Robinson theories and
universal domains. The real work starts with Lemma 5.10 of \cite{CHP}, which out of a
group configuration, produces a quantifier-free definable subgroup of an
algebraic group, in some reduct $\Omega[m]$. Note that in the cumulative
case, the subgroup $G_1$ can be chosen so that its generic type is
cumulative, by Proposition 1.15 of \cite{CHld}.  \
Then all results of \cite{CHP} up to Proposition 5.14 go through without
change.\\[0.05in]
(5.15) is the statement of the trichotomy theorem:\\[0.05in]
{\bf Theorem 5.15}. {\em Let $p$ be a basic type, and assume that
  $\calx_p(K)$ is not modular. Then $\calx_p(K)$ interprets an
  algebraically closed field of rank $1$.}\\[0.05in] 
The proof given in \cite{CHP} goes through, as it is just an adaptation
of the   proof of \cite{HZ} to our particular case. \\[0.1in]
We now come to the main result of the paper, given at the beginning of
section 6:\\[0.05in]
{\bf Theorem}. {\em Let $K\models {\rm ACFA}$, let $E=\acls(E)\subseteq
  K$, and let $p$ be a type over $E$, with $\SU(p)=1$. Then $p$ is not
  modular if and only if $p$ is non-orthogonal to the formula
  $\si^m(x)=x^{p^n}$ for some relatively prime $m,n\in \zee$ with $m\neq 0$.}\\[0.05in]
The proof goes through verbatim, to show that for some $m>0$, (passing
maybe to a larger  $E$), if $a$ realises $p$, there is some
$a'\in\acls(Ea)$ such that $\evSU(a'/E)=\SU(a'/E)[m]=1$, and $\qftp(a'/E)[m]$ is non-orthogonal to the formula $(\si^m)^r(x)=\frob^n(x)$ for some
integers $r\neq 0$ and $n$, with $(n,r)=1$ (and in fact, $r=1$). The
proof is now routine, using Lemma 1.12 of \cite{CH}: let $b,c$ be tuples
such that, in $\Omega[m]$, $c$ is independent from $\acls(Ea)=\acls(Ea')$ over $E$, $b$ satisfies
$(\si^m)^r(x)=\frob^n(x)$ and belongs to
$E_0=\acl_{\si^m}(Ea'c)$. The proof of Lemma 1.12 of \cite{CH} then
gives us an
$\acls(Ea)$-$\si^m$-embedding $\varphi$ of $F_0=\acls(Ea)E_0$ into $\Omega[m]$, such that  the fields $\si^i\varphi(F_0)$,
$i=0,\ldots,m-1$ are linearly disjoint over $\acls(Ea)$. It then follows
that $\varphi(c)$ is independent from $a$ over $E$ (in $\Omega$), and therefore
$p$ is non-orthogonal to $\si^{mr}(x)=\frob^n(x)$. \\[0.05in]
The proofs of the results of chapter 7 are also unchanged. \\[0.2in]

\bigskip

We have proved one part of the trichotomy, namely the dichotomy between modularity and a field structure.  
The second leg is proved in all characteristics in \cite{CH}, 5.12:  if $p$ is modular but has nontrivial algebraic closure geometry,
 then $p$ is non-orthogonal to an SU-rank one
  definable subgroup of an algebraic group, indeed of the additive or
  multiplicative group, or a simple abelian variety.   

Additional information concerning the non-orthogonality is available, see \cite{CH2}.  The internal structure 
of modular subgroups of semi-abelian varieties is fully understood, \cite{CH3}.     In the additive case, a bilinear map is definable in some cases; describing the full induced
structure remains open.   
 
  %\red{The discussion done in paragraph (5.20) of  \cite{CHP} then gives the full trichotomy theorem: }
 
%
%\noindent{\bf Trichotomy Theorem} \red{\em Let $p$ be a type of SU-rank $1$. Then
%either }
%
%\red{\em Moreover, if $p$ is non-orthogonal to a group, it is
%  non-orthogonal to the generic of SU-rank $1$ definable subgroup of an algebraic group, which is either
%  ${\mathbb G}_a$, ${\mathbb G}_m$ or a simple abelian variety. For a
%direct proof, the reader can see \cite{HMM}, Chapter 4.} 

%% One can similarly prove that  if $p$ is locally modular and the  combinatorial geometry associated to algebraic closure on $p$ is not trivial, then $p$ is non-orthogonal to an SU-rank one
%% definable subgroup of an algebraic group, indeed of the additive or
%% multiplicative group, or a simple abelian variety.  For an analysis of locally modular types of this kind, see \cite{HMM} and  \cite{CH3}.

\noindent 
%% {\bf Addresses of the authors}\\[0.1in]
%% CNRS (UMR 8553) - Ecole Normale Sup\'erieure\\
%% 45 rue d'Ulm\\
%% 75230 Paris cedex 05\\
%% France\\[0.1in]
%% Department of Mathematics\\
%% Hebrew University\\
%% Jerusalem, Israel 91904
\end{document}